\documentstyle[aps,pre,preprint,lscape,epsfig,subfigure,multirow]{revtex}
\def\f{\frac}
\def\l{\left}
\def\r{\right}

\begin{document}
%\begin{titlepage}

\title{Shock-capturing with natural high frequency
oscillations}

\author{Y. C. Zhou, Yun Gu and G. W. Wei
\footnote{Corresponding author, cscweigw@nus.edu.sg} \\  
Department of Computational Science \\
National University of Singapore, Singapore 117543}
\date{\today}
\maketitle
\begin{abstract}

This paper explores the potential of a newly developed 
conjugate filter oscillation reduction (CFOR) scheme
for shock-capturing under the influence of natural 
high-frequency oscillations. The conjugate low-pass and 
high-pass filters are constructed based on the principle
of the discrete singular convolution. Two Euler systems, 
the advection of an isentropy vortex flow and the interaction
of shock-entropy wave are considered to demonstrate the utility
of the CFOR scheme. Computational accuracy and order of 
approximation are examined and compared with the 
literature. Some of the best numerical results are 
obtained for the  shock-entropy wave interaction. 
Numerical experiments indicate that the proposed scheme is
stable, conservative and reliable for the numerical simulation 
of hyperbolic conservation laws.
 
\end{abstract}

%\newpage

\section{Introduction}

Scientists and engineers often face challenging problems in 
scientific and engineering computing. Some of typical problems 
include numerically induced chaos due to initial value being 
close to homoclinic manifold singularity in nonlinear 
Hamiltonian systems. Numerical instability in the prediction 
of higher-order eigenvalues is a long standing obstacle for 
the progress in the optimization of space structures.
In aerodynamics and fluid dynamics, one of major difficulties 
is to attain solutions that are free of spurious oscillations for 
compressible Euler equations involving steep gradient at 
discontinuities. In particular, to construct a high-accuracy 
and high-resolution solution for a system involving the 
interaction of shock-turbulent boundary layer is a severe 
challenge due to its natural high-frequency components. 
Traditional schemes such as upwind, Riemann solver, 
approximate Riemann solver, random choice method, 
artificial viscosity, are usually of low order in nature. 
More sophisticated approaches, such as total variation 
diminishing (TVD), essentially non-oscillatory 
(ENO)\cite{shu1}, weighted essentially non-oscillatory (WENO),
characteristic-based-split (CBS)\cite{zienkwicz99}
and discontinuous Galerkin schemes\cite{Galerkin2}
are proposed in the past two decades. Generally, these 
shock-capturing schemes are shown to be very successful 
in many applications.

Recently, it has been pointed out\cite{LLM} that the use of a sixth 
order accurate ENO scheme in the entire computational domain 
leads to a significant damping of turbulent fluctuations.
Garnier et al\cite{GMSCD} found that in the framework of freely 
decaying turbulence, the numerical dissipation of high-order 
accurate shock-capturing schemes masks the effect of the 
subgrid-scale model. Therefore, alternative approaches 
are of pressing desirable for many practical applications. 
The use of post-processing based filtering is one of 
important alternative approaches which can overcome the problem
of excessive numerical dissipation in many sophisticated 
shock-capturing schemes. Engquist et al\cite{engquist} 
proposed a set of nonlinear filters which discriminate and 
eliminate the dispersive wiggles in the basic solution.
Recently, Garnier et al\cite{eric01} have reported the use 
of the nonlinear dissipation components in some high-order 
shock-capturing schemes, such as the ENO and WENO, as filters.
With appropriate sensors, these ENO filters are shown to 
effectively improve the resolution of density waves,
entropy waves and stochastic turbulent fluctuations. 
Indeed, in the case of direct simulation of turbulence, 
typical successful numerical approaches are spectral 
methods for their capability of resolving multiscale 
features in turbulent fluctuations. However, spectral 
methods are notorious for their Gibbs' oscillations at 
the discontinuity. Therefore, it is highly desirable
to have methods that are of (arbitrary) high accuracy for 
resolving high-frequency waves and stochastic turbulent 
fluctuations, and are capable of shock-capturing without 
excessive numerical dissipation.

Conjugate filter oscillation reduction (CFOR)\cite{weigu,guwei} 
is one of such schemes newly developed for solving practical 
problems. The CFOR scheme is constructed based on the discrete 
singular convolution (DSC) algorithm\cite{weijcp,weijpa,weicmame}, 
a new approach for the numerical computation of singular 
convolutions. The theoretical foundation of the DSC  algorithm
is the theory of distributions and the theory of wavelet analyses. 
The DSC algorithm provides a unified approach to conventional 
local and global methods and has controllable accuracy for 
numerical solution of differential equations. The essential idea 
in the CFOR scheme is to use  conjugate (DSC) low-pass filters 
to remove the spurious oscillation generated by conjugate 
(DSC) high-pass filters which are implemented for the 
numerical approximation of differentiation operators.  
Conjugate filters are constructed by using DSC kernels 
and are optimal in the sense that they have similar order 
of regularity, accuracy, frequency bandwidth and computational 
supports. The CFOR scheme has been successfully applied to 
shock-capturing in association with Burgers' equation, 
one- and two-dimensional (2D) Euler systems including 
the Sod and Lax problems, and the Mach 3 flow past a wind 
tunnel with a step. The most promising feature of the 
CFOR scheme is that, the approach has controllable order of 
approximation for shock-capturing under practical situations. 
The objective of the present work is to explore the utility 
and limitation of the CFOR scheme in dealing with the
problem of shock-high frequency entropy wave interaction,
which is very challenging because conventional methods encounter
the difficulty of either insufficient accuracy or excessive 
numerical damping. It is believed that a better understanding
of the CFOR scheme is of importance to the development of 
high-accuracy and low-dissipation schemes for the numerical 
solution of more challenging problems, such as shock-turbulence 
interaction.

This paper is organized as follows. A brief retrospection is 
given to the DSC algorithm and the CFOR scheme in Section II.
Section III is devoted to the numerical experiments of 
1D and 2D Euler systems, an isentropic vortex flow and the
interaction of  shock-entropy wave. The latter is designed 
to test the capability of resolving shock from high-frequency 
entropy waves. A conclusion ends the paper.

\section{The DSC algorithm and CFOR scheme}

\subsection{DSC filters}

Singular convolutions occur commonly in science and engineering. 
Discrete singular convolution (DSC) is an effective  approach 
for the numerical realization of singular convolutions. 
There are many detailed descriptions about the discrete 
singular convolution in the literature\cite{weijcp,weijpa,weicmame}. 
The introduction in Ref. \cite{weijcp} is recommended for 
its theoretical underpinning and approximation philosophy. 
For the sake of integrity and convenience, a brief review of 
the DSC algorithm is given before describing the CFOR scheme.

In the context of distribution theory, a singular 
convolution can be defined by
\begin{eqnarray}
F(t)=(T * \eta)(t) = \int_{-\infty}^{\infty} T(t-x) \eta (x)dx, 
\label{dsc1} 
\end{eqnarray}
where $T$ is a singular kernel and  $\eta (x)$ is an element of 
the space of test functions. Interesting examples include singular
kernels of Hilbert (and Abel) type and delta type. The former 
plays an important role in the theory of analytical functions, 
processing of analytical signals, theory of linear responses 
and Radon transform. Since delta type kernels are the 
key element in the theory of approximation and the numerical 
solution of differential equations, we focus on the singular 
kernels of delta type
\begin{eqnarray}
T(x)= \delta^{(q)}(x),  ~~~~ (q=0, 1, 2,\cdots),  
\label{deltakn}
\end{eqnarray}
where superscript $(q)$ denotes the $q$th-order ``derivative'' 
of the delta distribution, $\delta (x)$, with respect to $x$, 
which should be understood as generalized derivatives of 
distributions. When $q=0$, the kernel, $T(x)=\delta(x)$, 
is important for the interpolation of surfaces and curves, including
applications to the design of engineering structures. For  hyperbolic
conservation law and Euler systems, two special cases, $q=0$ and $q=1$
are involved, whereas for the full Navier-Stokes equations, 
the case of $q=2$ will be also invoked.  Because of its singular
nature, the singular convolution of Eq. (\ref{dsc1}) cannot 
be directly used for numerical computations. In addition, the 
restriction to the test function is too strict for most practical 
applications. To avoid the difficulty of using singular 
expressions directly in numerical computations, we consider a 
discrete singular convolution which provides appropriate 
approximations to the original distribution
\begin{eqnarray}
f^{(q)} (x) \approx \sum_{k} \delta^{(q)}_{\Delta}(x-x_{k}) f(x_{k}),
 \label{fdsc1}
\end{eqnarray}
where $ \delta^{(q)}_{\Delta}(x-x_{k})$ are approximations 
to $\delta^{(q)}(x-x_{k})$ and are designed for being used 
in (discrete) summations. Here, $\{ x_{k} \}$ is an appropriate 
set of discrete points centered around the point $x$ and $\Delta$ 
is the grid spacing. Depending on the mathematical properties 
of the kernel, $\delta_{\Delta}$, the restriction on the 
$f(x_k)$ can be relaxed to include many common-occurring functions.
A variety of candidates are available for  $\delta_{\Delta}$
in the literature. In these examples, Shannon's delta kernel 
is of particularly interesting
\begin{eqnarray}
\delta_{\Delta}(x)=\text{sin}(\alpha x)/\pi x.    
\label{sdsc}
\end{eqnarray}
Shannon's kernel is a delta sequence and thus provides 
an approximation to the delta distribution
\begin{eqnarray}
\lim_{\alpha \rightarrow 0} < \f{\text{sin}(\pi x/\Delta)}{\pi x}, 
\eta (x) > = \eta(0). 
\label{deltas}
\end{eqnarray}
Shannon's kernel has been widely used in information theory, 
signal and image processing because the Fourier transform of 
Shannon's kernel is an ideal low-pass filter. However, 
the use of Shannon's kernel is limited by the fact that it 
has a slow-decaying oscillatory tail proportional to $\f{1}{x}$ 
in the coordinate domain. For signal processing, 
Shannon's kernel is an infinite impulse response (IIR) low-pass
filter. Therefore, when truncated in computational applications, 
its Fourier transform contains evident oscillations. A cure to 
this problem is to regularize Shannon's kernel with a Gaussian
\begin{eqnarray}
\delta_{\sigma, \Delta}(x) = 
    \f{\text{sin}(\pi x/\Delta)}{\pi x}
    e^{-\f{x^2}{2\sigma^2}}, \quad \sigma >0. 
                             \label{rgshannon}
\end{eqnarray}
Since $e^{-\f{x^2}{2\sigma^2}}$ is a Schwarz class function, 
it makes the regularized Shannon's kernel applicable 
to tempered distributions. Moreover, as the regularized kernels 
decay very fast in the space domain, they can be utilized as 
finite impulse response (FIR) low-pass filters. Their 
oscillation in the Fourier domain is dramatically reduced
and effectively controlled.

For sequences of the delta type, an interpolating algorithm 
sampling at the Nyquist frequency, $r = \pi/\Delta$, has 
an advantage over a non-interpolating discretization. 
Therefore, on a uniform grid, the regularized Shannon's 
kernel is discretized as
\begin{eqnarray}
\delta_{\sigma,\Delta}(x-x_k)=\f {\sin \f {\pi}{\Delta} (x-x_k)}
{ \f{\pi}{\Delta}(x-x_k)} e^{- \f {(x-x_k)^{2}}{2 \sigma^{2}} }
 \label{sdelta}.
\end{eqnarray}
The regularized kernel $\delta_{\sigma,\Delta}(x)$ corresponds 
to a family of FIR low-pass filters, each with a different 
compact support, according to $\f{\sigma}{\Delta}$, in the 
coordinate domain. Its $q$th order derivative is given 
by analytical differentiation
\begin{eqnarray}
\delta^{(q)}_{\sigma,\Delta}(x-x_k)=\l( \f{d}{dx} \r)^{(q)}
\f {\sin \f {\pi}{\Delta} (x-x_k)}
{ \f{\pi}{\Delta}(x-x_k)} e^{- \f {(x-x_k)^{2}}{2 \sigma^{2}} }.
 \label{sdeltadiff}
\end{eqnarray}
In this work, $\delta^{(q)}_{\sigma, \Delta}(x),~(q=0,1, \cdots)$ 
are referred as a family of ``conjugate filters", as they 
are derived from one generating function and consequently 
have similar degree of regularity, smoothness, time-frequency 
localization, effective support and bandwidth.

In application, optimal results are usually obtained if the 
window size $\sigma$ varies as a function of the central frequency 
$\pi / \Delta$, such that $r=\sigma / \Delta$ is a parameter 
chosen in computations. Both interpolation and differentiation 
are realized by the following convolution algorithm
\begin{eqnarray}
f^{(q)}(x) \approx \sum_{k=-W}^{k=W} \delta^{(q)}_{\sigma,\Delta}(x-x_k)f(x_k),
              \quad (q=0,1,2,\cdots),
\end{eqnarray}              
where $2W+1$ is the computational bandwidth, or effective kernel 
support, which is usually smaller than the entire computational 
domain, [a,b]. Expressions of $\delta^{(q)}_{\sigma,\Delta}(x)$ 
with $q=0,1$ are given as
\begin{eqnarray}
\delta_{\sigma,\Delta} (x) =
\l \{  
\begin{array}{ll}
  \f {\sin \l( \f {\pi x}{\Delta}   \r)  
      \exp \l(-\f {x^2}{2 \sigma^2} \r)} 
      {\f {\pi x}{\Delta}}      &   ~~~  (x \ne 0) \\ 
       1                        &   ~~~  (x=0), 
\end{array}
\r.                             \label{deta0}
\end{eqnarray}     
\begin{eqnarray}
\delta_{\sigma,\Delta}^{(1)} (x) =
\l \{  
   \begin{array}{ll}
    \f {\cos \l( \f  {\pi x}{\Delta} \r) 
        \exp \l( -\f {x^2}{2 \sigma^2} \r)}  {x} 
  - \f {\sin \l( \f  {\pi x}{\Delta} \r)
        \exp \l( -\f {x^2}{2 \sigma^2} \r)}  {\f {\pi x^{2}}{\Delta}}
          \\ 
  - \f {\sin \l( \f  {\pi x}{\Delta} \r) 
        \exp \l( -\f {x^2}{2 \sigma^2} \r)}  {\f {\pi \sigma^{2}}{\Delta}}
                     &   ~~~           (x \ne 0)   \\ 
   0                 &   ~~~         (x=0).
\end{array}
\r.                              \label{deta1}
\end{eqnarray} 
Expressions of higher-order derivatives 
for $\delta_{\sigma,\Delta}(x)$ can be found 
elsewhere\cite{weicmame}.

\subsection{The CFOR scheme}

Consider  2D Euler equations for gas dynamics 
in a vector notation having the conservation form of 
\begin{eqnarray}
U_t + F(U)_x + G(U)_y=0 \label{Euler}
\end{eqnarray}
with
\begin{eqnarray}
U= \l ( \begin{array}{c}
         \rho   \\
         \rho u \\
         \rho v \\
              E \\
       \end{array} \r ); \quad
F(U) = \l ( \begin{array}{c}
             \rho u       \\
             \rho u^2 + p \\
             \rho u v     \\ 
             u(E+p)  \\
           \end{array} \r ); \quad
G(U) = \l ( \begin{array}{c}
             \rho v       \\
             \rho u v     \\
             \rho v^2 + p \\ 
             v (E+p)  \\
           \end{array} \r ),          
\end{eqnarray}
where, $\rho, u, v, p$ and $ E$ denote the density, the velocities 
in  $x$- and $y$-directions, the pressure and the total energy per 
unit mass $E=\rho(e+(u^2+v^2)/2)$, respectively. Here, $e$ is the 
specific internal energy. For an ideal gas with constant 
specific heat ratio ($\gamma =1.4$) considered here, one has
$e=\f{p}{(\gamma-1)\rho}$.

Let denote spatial discretizations of 
$F(U)$ and $G(U)$ at a grid point $(i,j)$, 
as $F(U_{i,j})$ and $G(U_{i,j})$.
Their spatial derivatives $F(U)_x$ and $G(U)_y$
are approximated by the DSC high-pass filters 
according to Eq. (\ref{fdsc1}), i.e.
\begin{eqnarray}
F(U_{i,j})_x = 
\sum_{k=i-W}^{i+W} \delta_{\sigma,\Delta}^{(1)}(x_i-x_k) F(U_{k,j})
    \label{spdisrx} \\
G(U_{i,j})_y = 
\sum_{k=j-W}^{j+W} \delta_{\sigma,\Delta}^{(1)}(y_j-y_k) G(U_{i,k}).
    \label{spdisry}
\end{eqnarray}
The accuracy of the DSC algorithm is controllable\cite{weicmame}.

Although there is no rigorous proof about whether the 
standard forth-order Runge-Kutta (RK-4) scheme is TVD, it is 
still one of the most widely used temporal scheme for 
hyperbolic conservation laws. The RK-4 is adopted in the present 
work and it takes the following form in our problems
\begin{eqnarray}
k_1 & = & - F(U^n_{i,j})_x - G(U^n_{i,j})_y        \\
k_2 & = & - F(U^n_{i,j}+\f{\Delta t}{2} k_1)_x 
          - G(U^n_{i,j}+\f{\Delta t}{2} k_1)_y \\
k_3 & = & - F(U^n_{i,j}+\f{\Delta t}{2} k_2)_x 
          - G(U^n_{i,j}+\f{\Delta t}{2} k_2)_y \\
k_4 & = & - F(U^n_{i,j}+\Delta t k_3)_x 
          - G(U^n_{i,j}+\Delta t k_3)_y      \\
U^{n+1}_{i,j} & = & U^n_{i,j} + \f{\Delta t}{6} [k_1+2k_2+2k_3+k_4]. 
\label{basic}
\end{eqnarray}
The approximation of  $F(U)$ and $G(U)$ by using the DSC high-pass 
filters in Eq. (\ref{spdisry}), together with the RK-4
scheme (\ref{basic}), provides a basic scheme for the 
numerical integration of the Euler system, Eq. (\ref{Euler}).  
No additional effort is required if the problem under 
consideration does not involve discontinuity.  
Otherwise, an additional filtering can be implemented to 
prevent spurious oscillations.

As aforementioned, the conjugate filters are constructed 
from the same generating function (\ref{sdelta}). The 
frequency responses of conjugate low-pass  and high-pass 
filters  are illustrated  in Fig. \ref{dscdft}. It can be 
seen that below $0.7\pi / \Delta$, all the conjugate filters 
are highly accurate. However, in the high-frequency region, 
the frequency responses of both the low-pass filter and 
first-order high-pass filter are serious under estimating, 
whereas the frequency response of the second-order high-pass 
filter is over estimating. The error in the  high-frequency
response is harmless for numerical problems involving only 
low frequency components. However, in the case of shock 
and discontinuity, the solution contains much high-frequency 
component, the error in the high-frequency response will be 
accumulated and amplified during the time integration, and 
leads to spurious oscillations. This observation motivates 
us to use the conjugate low-pass filter to appropriately 
eliminate most of the high frequency response produced by 
the conjugate high-pass filters. As a result, the solution 
generated by conjugate filters is reliable for the 
frequency below the effective bandwidth of the filters.
The effective bandwidth or frequency cut-off is controlled 
by the choice of the DSC parameter $r=\sigma/\Delta$, for a 
given $\Delta$. The CFOR scheme is implemented via the 
following two-step procedure
\begin{eqnarray}
\hat U^{n+1}_{i,j}  & = & H(U^n_{i,j}) \\
U^{n+1}_{i,j} & = & L(\hat U^{n+1}_{i,j}),
\end{eqnarray}
where $H(U^n_{i,j})$ is the high-pass filtering process
given by $ U^n_{i,j} + \f{\Delta t}{6} [k_1+2k_2+2k_3+k_4]$,
i.e., the intermediate result obtained by using the scheme 
(\ref{basic}). Here, $L$ is the DSC low-pass filtering as 
shown in Eq. (\ref{fdsc1}) with $q=0$. This interpolative 
low-pass filter is implemented through prediction (in 
which the variables on the grid are interpolated to 
the middle points of the cells) and restoration (in which 
the variables on the grid are restored from their values 
at the middle points of the cells)\cite{guwei}. 

In the above two-step procedure, the second step, i.e. the 
low-pass conjugate filtering is controlled (turned on or 
turned off) by a sensor. There are a number of such sensors 
which can be used in the present scheme. Among them, the TVD 
switch seems to be the simplest one. In this work, we defined 
a high-frequency measure $W$. Upon the increment $\delta W$ 
exceeding a prescribed threshold $\zeta$, the low-pass filter process 
is carried out and $U^{n+1}_{i,j}$ is the solution in the new 
time step $n+1$. Otherwise, no low-pass filter will be exerted 
to $\hat U^{n+1}$ and the latter is taken as the result at the
$n+1$ time step. The high-frequency  measure $W$ is defined via 
a multiscale wavelet transform of a set of discrete function 
values at time $t_n$ as 
\begin{eqnarray}
||W^{n}||=\sum_{m}||W^n_m||, 
\end{eqnarray}
where $||W^n_m||$ is given by a convolution with a 
wavelet $\psi_{mj}$ of scale $m$
\begin{eqnarray}
||W^n_m||=\sum_{k}|\sum_{j} \psi_{mj}(x_k)u^n(x_j)|.
\end{eqnarray}
Such a definition can be further illustrated by one of its 
special case - the TVD sensor, which can be obtained by 
restricting to the Haar wavelet with a single scale. 
The choice of the threshold $\zeta$ for the 
high-frequency measure $\delta W = || W^{n+1} ||- || W^n|| $ 
depends on the nature of the problem under study and is in 
the range from $0.001$ to $0.002$.

\section{Results and discussions}

In this section, we examine the utility and explore 
limitation of the proposed CFOR scheme by using two
benchmark numerical problems, the 2D advection 
of an isentropic vortex\cite{shu1,eric01} and the 
interaction of shock-entropy wave\cite{shu1}. 
The first example is designed to quantitatively 
access the phase and amplitude errors of the CFOR
scheme in handling 2D Euler problems. To maintain a 
small error in both phase and amplitude is particularly 
desirable for a scheme to handle shock-entropy wave 
interaction and many other aerodynamic problems.
Moreover, extensive numerical data are available for this 
problem and a comparison with many other shock-capturing 
schemes, such as the ENO and WENO, is readily possible.  
The second problem is a standard test for the numerical 
ability of treating high-frequency entropy wave-shock
interaction. It is a severe challenge for most existing 
shock-capturing methods due to its high-frequency 
nature. Numerical results can be objectively 
evaluated by a quantitative criterion obtained from a 
linear analysis. Parameters $W=32$ and $r=3.2$ are 
used for all the high-pass filtering and low-pass 
prediction. For the low-pass restoration, $r=3.2$ 
is used in the evolution of isentropic vortex, while 
$r$ values of $1.9\sim 2.1$ are used in the interaction 
of shock and high-frequency entropy waves.

\subsection{Isentropic vortex}

To quantitatively analyze the performance of the proposed
CFOR scheme, the advection of an isentropy vortex in a free
stream is computed. As the exact solution of the problem 
is available, it is an excellent benchmark for accessing 
the accuracy and stability of shock-capturing schemes  
and has been previous considered by many 
researchers\cite{shu1,eric01}.

Consider a mean flow of $(\rho_{\infty},u_{\infty},v_{\infty},
P_{\infty},T_{\infty})=(1,1,1,1,1)$ with a periodic boundary 
condition in both directions. At $t_0$, the flow is 
perturbed by an isentropic vortex $(u',v',T')$ 
centered at $(x_0,y_0)$, having the form of 
\begin{eqnarray}
u' & = & -\f{\lambda}{2 \pi}(y-y_0) e^{\eta(1-r^2)}, \\
v' & = & ~~\f{\lambda}{2 \pi}(x-x_0) e^{\eta(1-r^2)}, \\
T' & = & -\f{(\gamma-1) \lambda^2}{16 \eta \gamma \pi^2} e^{2 \eta (1-r^2)}.
\end{eqnarray}
Here, $r=\sqrt{(x-x_0)^2+(y-y_0)^2}$ is the distance to the 
vortex center; $\lambda$ is the strength of the vortex and 
$\eta$ is a parameter determining the gradient of the solution, 
and is unity in this study. Note that for an isentropic flow, 
relations $p=\rho^{\gamma}$ and $T=p/\rho$ are valid. 
Therefore, the perturbation in $\rho$ is required to be
\begin{eqnarray}
\rho=(T_{\infty} + T')^{1/(\gamma-1)} =
\l [ 1-\f{(\gamma-1) \lambda^2}{16 \eta \gamma \pi^2} e^{2 \eta (1-r^2)}
\r ] ^{1/(\gamma-1)}. 
\label{density}
\end{eqnarray}
For a comparison with the existing literature\cite{eric01},
the computational domain is chosen as $[0,10]\times[0,10]$ 
with the center of the vortex being initially located at 
$(x_0,y_0)=(5,5)$, the geometrical center of the computational
domain. Two experiments are performed in this study. One is 
to examine the accuracy of the CFOR scheme and to compare with 
the available literature. The other is to investigate the 
stability and performance of the CFOR scheme for long-time 
integration. For the first experiment, we compute the 
density profile up to $t=2$ using five sets of meshes 
($N=20,40,80,160,320$) which are selected by Garnier et 
al\cite{eric01}. In the present computations, the CFL number 
is chosen as $0.5$ for a comparison with previous 
results\cite{eric01}.

Two error measures, $L_1$ and $L_2$, are used in this
study. To be consistent with the literature\cite{eric01},
two errors used in this paper are defined as
\begin{eqnarray}
L_1 & = & {1\over (N+1)^2} 
  \sum_{i=0}^{N} \sum_{j=0}^{N} | f_{i,j}- \bar{f}_{i,j} |  \\
L_2 & = & {1\over (N+1)} \sqrt{\sum_{i=0}^{N} \sum_{j=0}^{N}
               | f_{i,j}- \bar{f}_{i,j} |^2},
\end{eqnarray}
where $f$ is the numerical result and $\bar{f}$ the exact 
solution (Note that they are not the standard definitions).
The CFOR errors for the density with respect to the 
exact solution are listed in Tables \ref{vortextable1} and 
\ref{vortextable2}. Highly accurate results are obtained, 
as shown by the tables. Obviously, the proposed scheme is 
much more accurate than any other scheme listed in the tables, 
which are reported  by Garnier et al\cite{eric01}.
Remarkably, the CFOR scheme is from 4 to 5 orders more 
accurate than other schemes when $N=80$.

As the spatial discretization of the CFOR scheme
is extremely accurate, some of the present results
computed at CFL=0.5 might be limited by the CFL 
number, which was optimized according to various schemes 
given in Ref. \cite{eric01}. This is indeed the case. 
The CFOR results computed at CFL=0.01 are generally 
more accurate than those obtained at  CFL=0.5 as shown 
in Tables \ref{vortextable1} and \ref{vortextable2}. Note 
that the accuracy of the CFOR scheme increases dramatically 
when the mesh is refined from 40 grid points to 80 grid 
points, with the numerically computed approximation order 
being more than 15. Therefore, the CFOR scheme has the feature 
of spectral-like methods. Obvious, due to its extremely
high accuracy, the CFOR can be used for large scale 
simulations without resorting to a very large mesh as 
required by low order schemes.

Our second numerical experiment concerns the performance
of the CFOR scheme for the long-time integration, which 
poses a severe challenge to the stability and conservation 
of the discretization scheme\cite{shu1}. The solution of 
$\rho$ is sampled at $t=2, 10, 50$ and $100$, with the grid 
spacing of $\Delta x= \Delta y=0.125$ and  CFL=0.5. 
In Fig. \ref{vortexfig1}, we show the horizontal line cut 
through the center of the vortex for density $\rho$.
Obviously, there is no visual deviation between the 
computed result and exact result. Errors listed in Table 
\ref{vortextable3} further confirm that the present  
scheme is extreme accuracy, free of excessive dissipation
and reliable.

Since there is no presence of shock in this case, the low-pass 
filter originally designed to suppress dispersive wiggles 
might appear useless. In this experiment, it is found that 
the DSC algorithm on its own can already provide excellent 
results if the integration time is small enough. Thus,
the conjugate low-pass filter does not need to be activated 
during an initial time period. However, as the time progresses, 
errors would accumulate rapidly and the computation  could 
become unstable if the low-pass filter were not used to 
effectively control the dramatical nonlinear growth of the errors. 
Therefore, the CFOR scheme is very robust for the treatment 
of this problem. As shown in Table \ref{vortextable3} and 
Fig. \ref{vortexfig1} ($t=100$), the long time simulation 
results are very stable. The vortex core is well conserved 
and the accuracy is extremely high. These results indicate 
that the CFOR  scheme is highly accurate, stable and 
conservative for the long-time integration of Euler systems.
It is a potential approach for the numerical integration of 
hyperbolic conservation laws. Its ability for shock-capturing 
is examined in the next subsection.

\subsection{Interaction of shock and high-frequency entropy wave}

The interaction of shock and high frequency entropy wave
is a standard test problem for benchmarking potential high-order 
shock-capturing methods. The problem is significant due to
its relevant to the interaction of shock-turbulence. A 
Mach 3 right-moving shock interacts with a small amplitude
entropy wave. The computation domain is taken as $[0,5]$
and the flow field is initialized with
\begin{eqnarray} \label{inicond}
 (\rho, u, p)  = \left \{
         \begin{array}{lccl}
         (3.85714,                       & 2.629369, & 10.33333 ); &
                          x \le 0.5 \\
         (e^{- \epsilon \sin( \kappa x)}, & 0,        & 1.0 ); &
                          x > 0.5 \\
        \end{array} \right.,
\end{eqnarray}        
where $\epsilon$ and $\kappa$ are the amplitude and the wave 
number of the entropy wave before the shock. The amplitude and 
the wave number of amplified wave after the shock can be 
obtained from a linear analysis\cite{mckenzie}. In our 
numerical experiments, we vary the wave number of the 
pre-shock entropy wave while keep its amplitude unchanged.
As a results, the amplitude of the post-shock entropy wave 
will also be a constant, i.e. 0.08690716 and the corresponding 
amplitude of pre-shock entropy wave is $\epsilon$=0.01.

In this problem, a large-amplitude high-frequency entropy 
wave is mixed with spurious oscillations. It is difficult 
to distinguish them clearly in numerical simulations. 
Potential methods designed for suppressing the spurious 
oscillation might also smear the high-frequency post-shock 
entropy wave. As the wave number $\kappa$ increases, the problem 
becomes extremely challenging\cite{shu1}. Low-order 
shock-capturing schemes and even some  popular high-order 
schemes encounter the difficulty in preserving the amplitude of 
the entropy wave due to excessive dissipation with a given mesh size. 
Therefore, a success shock capturing method should be able to 
eliminate Gibbs' oscillation, capture the shock
and  preserve the entropy wave.

The computational domain is deployed with 800 grid points, and 
such a mesh is used in all the numerical tests except for further 
specified. First, we consider the case of $\kappa$=13. This 
is a good test case for a basic scheme as there are 20 grid 
points per entropy wavelength, which is sufficient for describing 
the wave if there is no shock. Such a case was found being 
slightly difficult for the fifth order WENO scheme\cite{shu1}. A 
significant amplitude damping occurs and the mesh of $N=1200$ has 
to be used to maintain the amplitude of the entropy wave\cite{shu1}. 
The result of the CFOR scheme is depicted in Fig. \ref{shockfig1}(a). 
It is  seen that the generated entropy wave spans fully over the 
strip bounded by two solid lines, showing excellent agreement 
with the linear analysis. The shock is exactly captured with a 
small frequency mismatch locates at the shock front. Such 
mismatch also occurs to the ENO and WENO schemes\cite{shu1} due 
to nature of discontinuity. The performance of the CFOR 
scheme is really remarkable for this case.

Next, we double the wave number, i.e. $\kappa$=26. The number 
of supporting grid points per generated entropy wavelength is 10. 
It is a quite difficult case for low-order shock-capturing 
schemes and no available result is reported in the literature,
to our knowledge.  The CFOR result is plotted in 
Fig. \ref{shockfig1}(b). Obviously, the compressed entropy wave 
is excellently resolved in the  post-shock regime. There is 
no visible trace of excessive dissipation as the amplitude of 
the entropy wave reaches its full strength in the whole 
post-shock regime. As expected, the frequency mismatch near 
the shock front becomes more obvious because there is an 
enlarged difference in the wave frequencies before and 
after the shock.

We further increase the wave number $\kappa$ to 39
and plot the CFOR result in  Fig. \ref{shockfig1}(c).
It is interesting to note that the compressed entropy 
wave peaks span to it full amplitude for only about half 
of its extrema over the post-shock regime. Analysis 
indicates that the presence of the under developed 
peaks is not due to excessive dissipation. Instead, it 
is due to the insufficient resolution in the plot. With a 
total of 800 grid points in the domain, there is less 
than 7 grid points per wavelength. Such a grid is not large
enough to fully resolve all the extrema in the compressed 
entropy wave. This explains the suppressed extrema in the 
plot. Obviously, it is extremely difficult to capture
shock on such a grid for any potential scheme. However, 
the CFOR scheme performs extremely well as shown in 
Fig. \ref{shockfig1}(d), which is obtained by interpolating 
the CFOR result in Fig. \ref{shockfig1}(c) to a denser grid 
$(N=1600)$. The interpolation is carried out by using the DSC 
interpolation scheme, i.e., the conjugate low-pass filter as 
given in Eq. (\ref{fdsc1}) with $q=0$. Apparently, the quality 
of this result is comparable with the case of $\kappa$=26. This 
confirms that the CFOR scheme works well for shock-capturing under 
a very small ratio of grid points and wavelength.

Finally,  we consider two large  wave numbers,
$\kappa=52$ and 65, to further test the performance of 
the CFOR scheme. A mesh of 800 grid points means 
a ratio of less than 5 grid points per wavelength, 
which is too few for simultaneous shock capturing and
high-frequency wave resolving. The CFOR scheme generates 
some small amplitude damping for $\kappa=52$ (which is 
not shown). Therefore, we increase the mesh size to 
$N=1600$ for these computations. With this mesh, there are 
10 grid points in each generated wavelength for $\kappa=52$. 
The resolution in the post-shock regime is excellent as shown 
in Fig. \ref{shockfig1}(e). As shown in  Fig. \ref{shockfig1}(f),
results for $\kappa=65$ are also very good.
However, there is a visible amplitude 
damping in the generated entropy waves. We noticed that the 
conjugate low-pass filter is activated more often in this 
computation than in previous cases. With the increase of 
wave number $\kappa$, the $r$ should also be increased to 
broaden the effective bandwidth so that high-frequency wave 
shall not suffer from too much numerical damping.
Obviously, the present scheme is capable of distinguishing 
the spurious oscillations from the high frequency oscillation
of the entropy wave, so as to attain high resolution for 
normal physical properties and capture the shock. 
It is believed that the potential of proposed CFOR scheme 
has been sufficiently demonstrated by these experiments.

\section{Conclusion}
 
The potential utility of a newly developed conjugate filter 
oscillation reduction (CFOR) scheme\cite{weigu,guwei} for 
the treatment of shock and high-frequency wave interaction 
is explored. The CFOR scheme is constructed based on the 
discrete singular convolution (DSC) algorithm, which is
a practical approach for the numerical realization of singular
convolutions. The essential idea of the CFOR scheme is to
employ a conjugate low-pass filter to effectively remove the 
high-frequency errors (spurious oscillations) created by
a set of high-pass filters, which are employed to discretize
the spatial derivatives in the hyperbolic conservation laws
or partial differential equations. The conjugate
low-pass and high-pass filters are optimal for shock-capturing
and spurious oscillation suppressing in the sense that they 
are generated from the same expression and consequently have 
similar order of regularity and approximation, effective 
frequency band and compact support.

The separation of the basic spatial discretization (the high-pass
filtering) and the post-processing (low-pass filtering) in the 
proposed shock-capturing scheme makes it possible to focus on
the design of a set of versatile and efficient filters. Such an 
approach has a few advantages. First, the basic DSC algorithm 
is a local method but it can be as accurate as a spectral method. 
Therefore, the CFOR scheme has controllable accuracy via the 
choice of filter parameters. Secondly, the post-processing applies 
at most once per iteration circle, comparing with pointwise 
treatment at each grid point in many other schemes. Therefore, 
there is a potential increase in the computational efficiency 
by a mature CFOR code. Finally, the implementation of the 
conjugate low-pass filter and filter parameters are easily 
controlled. Hence, physical high frequency can be effectively 
distinguished from the shock induced spurious oscillations 
in the present approach.

The performance of the proposed scheme is examined by using 
two benchmark Euler systems, the free evolution of a 2D 
isentropic vortex\cite{shu1,eric01} and the interaction of
shock-entropy wave\cite{shu1}. The first example has an exact
solution and its long-time evolution is a non-trivial
task. The CFOR scheme provides higher-order accuracy for 
solving the problem and its performance is compared with 
many other schemes in the literature\cite{shu1,eric01}.
The system of shock-entropy wave interaction is a difficult 
case due to its natural high frequency oscillations in 
the compressed entropy wave, which is easily damped by 
the excessive numerical dissipation in most existing 
shock-capturing schemes. The problem becomes a severe 
challenge as the wave frequency increases. It is demonstrated 
that the CFOR provides some of the best solution even 
available for this problem. The application of the CFOR 
scheme to more complicated problems and the adaptive 
optimization of the control parameters for conjugate 
filters are under consideration. 

\centerline{\bf Acknowledgment}

{This work was supported by the National University of 
Singapore. The authors thank Pierre Sagaut 
and Eric Garnier for useful corresponding.}

%\end{document}

\begin{center}
\begin{table}
\begin{tabular}{cccccccccccc}
N   &   & CFOR$^1$ & CFOR$^2$  & C4  & ENO  & MUSCL   
& WENO    & ENO$^{\it ACM}$  
& MUSCL$^{\it ACM}$ & WENO$^{\it ACM}$ \\ \hline
20  & error & 1.27E-3 & 1.27E-3 & 1.08E-2 & 7.83E-3 & 9.33E-3 & 6.12E-3 & 5.63E-3 & 6.18E-3 & 4.61E-3 \\ \hline
\multirow{2}{5mm}{40}  & error & 7.14E-6 & 3.76E-6 & 1.13E-3 & 1.28E-3 
      & 2.39E-3 & 9.39E-4 & 7.81E-4 & 1.29E-3 & 6.11E-4 \\ %\cline{2-11}
& order & 7.47 & 8.10 & 3.26 & 2.61 & 1.96 & 2.70 & 2.85 &  2.26 & 2.91 \\ \hline
\multirow{2}{5mm}{80}  & error & 4.57E-9 & 9.12E-11 & 5.78E-5 & 2.08E-4 
      & 5.99E-4 & 7.07E-5 & 6.68E-5 & 2.81E-4 & 4.58E-4 \\ %\cline{2-11}
& order & 10.61 & 15.33 & 4.29 & 2.62 & 1.99 & 3.73 & 3.55 & 2.19 & 3.74  \\ \hline
\multirow{2}{5mm}{160} & error & 3.23E-10 & 3.68E-11 & 3.79E-6 & 3.01E-5 
        & 1.26E-4 & 2.46E-6 & 7.84E-6 & 5.31E-5 & 2.95E-6 \\ %\cline{2-11}
& order & 3.82 & 1.31 & 3.93 & 2.79 & 2.25 & 4.84 & 3.09 & 2.40 & 3.97  \\ \hline
\multirow{2}{5mm}{320} & error & 5.09E-11 & 3.14E-11 & 2.41E-7 & 4.07E-6 
        & 2.26E-5 & 8.52E-8 & 6.82E-7 & 8.61E-6 & 2.13E-7 \\ %\cline{2-11}
& order & 2.67 & 0.23 & 3.97 & 2.89 & 2.47 & 4.85 & 3.52 & 2.62 & 3.79  \\
\end{tabular} 
\vspace{0.5cm}
\caption{
$L_1$ error for the density at t=2; 
The CFL number is 0.01 for CFOR$^2$ and
for all the other schemes, the CFL number is 0.5;
C4: fourth-order accurate, conservative centered scheme;
ENO: third-order;
MUSCL: third-order;
WENO: fifth-order;
XXX$^{\it ACM}$ denotes that the C4 scheme is used as the basic scheme 
with the XXX being 
a characteristic-based filter, while Harten's artificial compression 
method is used as a sensor to indicate the local numerical dissipation,
see Ref. [7]. }
\label{vortextable1}
\newpage

\begin{tabular}{cccccccccccc}
N   &  & CFOR$^1$ & CFOR$^2$ & C4  & ENO  & MUSCL   & WENO    & ENO$^{\it ACM}$  
& MUSCL$^{\it ACM}$ & WENO$^{\it ACM}$ \\ \hline
20  &  error & 3.28E-3 & 3.28E-3 & 1.93E-2 & 2.45E-2 & 2.90E-2 & 1.90E-2 & 1.77E-2 & 1.97E-2 & 1.45E-2 \\ \hline
\multirow{2}{5mm}{40}  &  error & 2.03E-5 & 1.74E-5 & 2.92E-3 & 4.09E-3 
      & 8.29E-3 & 3.16E-3 & 2.47E-3 & 4.05E-3 & 2.08E-3 \\  %\cline{2-11}
& order & 7.33 & 7.56 & 2.72 & 2.58 & 1.81 & 2.59 & 2.84 & 2.28 & 2.80 \\ \hline
\multirow{2}{5mm}{80}  & error & 1.47E-8 & 6.57E-10 & 1.90E-4 & 6.75E-4 
      & 2.26E-3 & 2.64E-4 & 2.08E-4 & 1.14E-3 & 1.48E-4 \\  %\cline{2-11}
& order & 10.43 & 14.69 & 3.94 & 2.60 & 1.88 & 3.58 & 3.57 & 1.83 & 3.81 \\ \hline
\multirow{2}{5mm}{160} & error & 1.05E-9 & 4.76E-10 & 1.23E-5 & 8.69E-5 
      & 5.91E-4 & 1.10E-5 & 2.51E-5 & 3.12E-4 & 9.44E-6 \\ %\cline{2-11}
& order & 3.81 & 0.46 & 3.95 & 2.96 & 1.94 & 4.58 & 3.05 & 1.87 & 3.97 \\ \hline
\multirow{2}{5mm}{320} & error & 4.17E-10 & 4.10E-10 & 7.84E-7 & 1.33E-5 
      & 1.31E-4 & 2.93E-7 & 2.19E-6 & 6.07E-5 & 6.85E-7 \\ %\cline{2-11}
& order & 1.33 & 0.21 & 3.97 & 2.71 & 2.17 & 5.23 & 3.52 & 2.36 & 3.78 \\
\end{tabular} 
\vspace{0.5cm}
\caption{ $L_2$ error for the density at t=2. See Table 
\ref{vortextable1} for  captions.} 
\label{vortextable2}
\vspace{3cm}

\begin{tabular}{ccccc}
Time  &  2  & 10  & 50  & 100              \\ \hline
$L_1$ &  4.56E-9 & 1.77E-8 & 4.27E-8 & 8.90E-8    \\ \hline
$L_2$ &  1.47E-8 & 4.95E-8 & 1.44E-7 & 3.01E-7    \\ 
\end{tabular} 
\vspace{0.5cm}
\caption{ Errors for the density at different times (CFL=0.5, N=80).}
\label{vortextable3}
\end{table}

\begin{figure}
 \centering  
 \includegraphics[width=8cm]{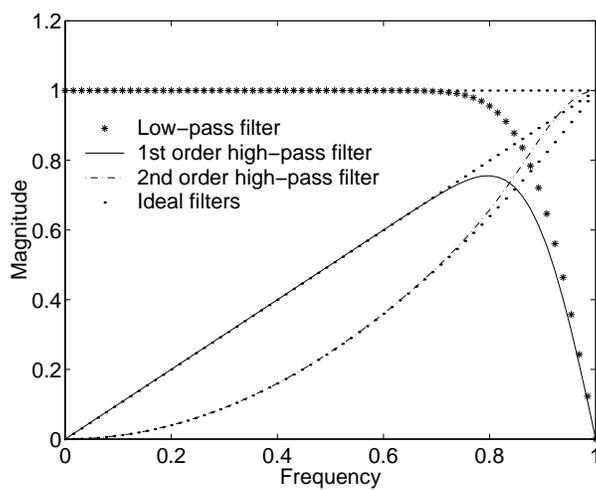}
 \caption{Frequency responses of the 
              conjugate DSC filters (in the unit of $\pi/\Delta$),
              The maximum amplitude of the filters is normalized 
              to the unit.} \label{dscdft}
\end{figure}

\newpage

\begin{figure}
 \centering  
 \subfigure[t=2]{
     \label{fig:subfig:a}
     \includegraphics[width=7cm]{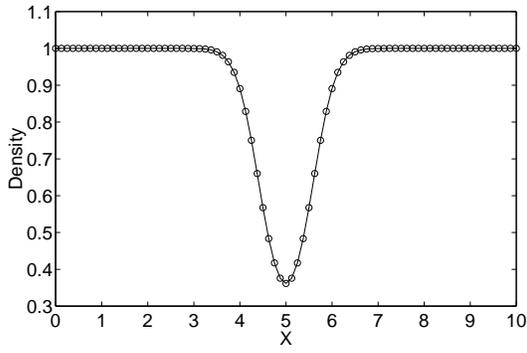}}
 \hspace{0.5cm} 
 \subfigure[t=10]{
     \label{fig:subfig:b}
     \includegraphics[width=7cm]{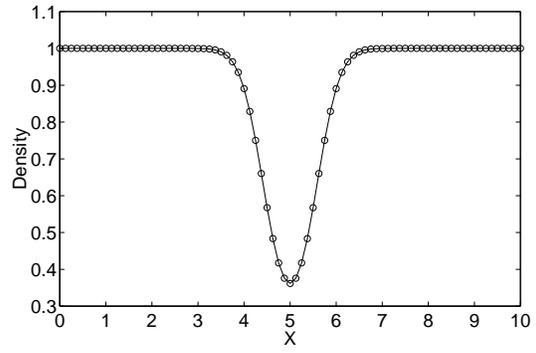}}  \\ %[5pt]
  \subfigure[t=50]{
     \label{fig:subfig:c}
     \includegraphics[width=7cm]{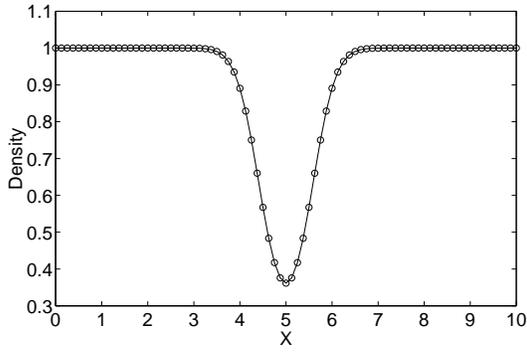}}
 \hspace{0.5cm}
 \subfigure[t=100]{
     \label{fig:subfig:d}
     \includegraphics[width=7cm]{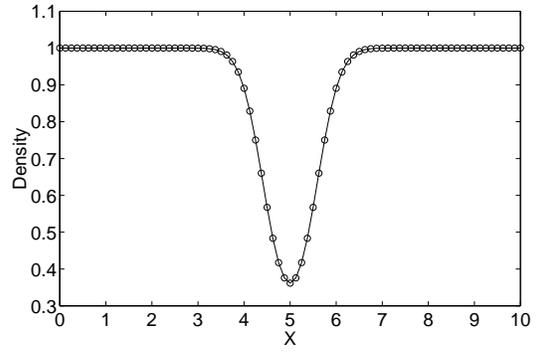}}
\caption{Density profiles in horizontal cutting at four times. 
Solid line is the exact profile and the circle dots 
denote the numerical results.}
\label{vortexfig1}
\end{figure}

\newpage

\begin{figure}
 \centering  
 \subfigure[$\kappa=13$, N=800]{
     \label{fig2:subfig:a}
     \includegraphics[width=12cm]{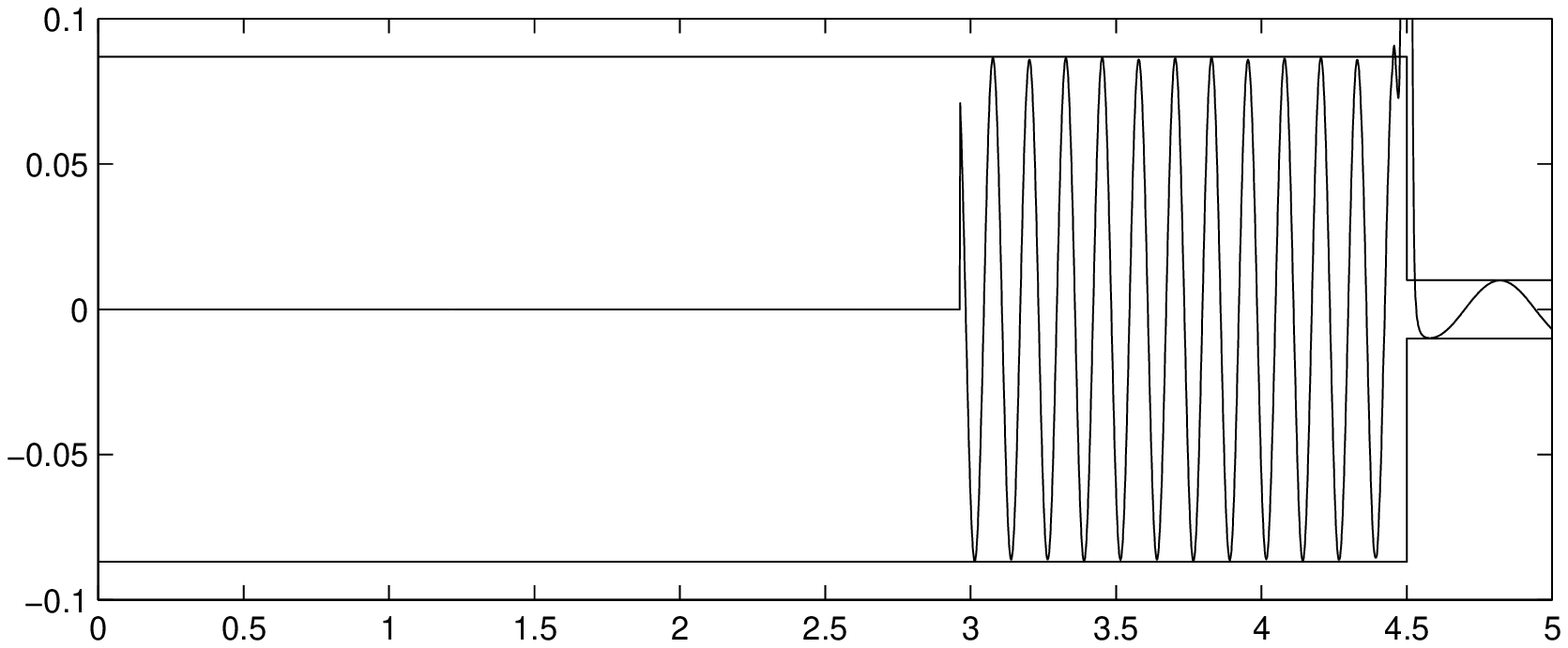}} \\ 
     
 \subfigure[$\kappa=26$, N=800]{
     \label{fig2:subfig:b}
     \includegraphics[width=12cm]{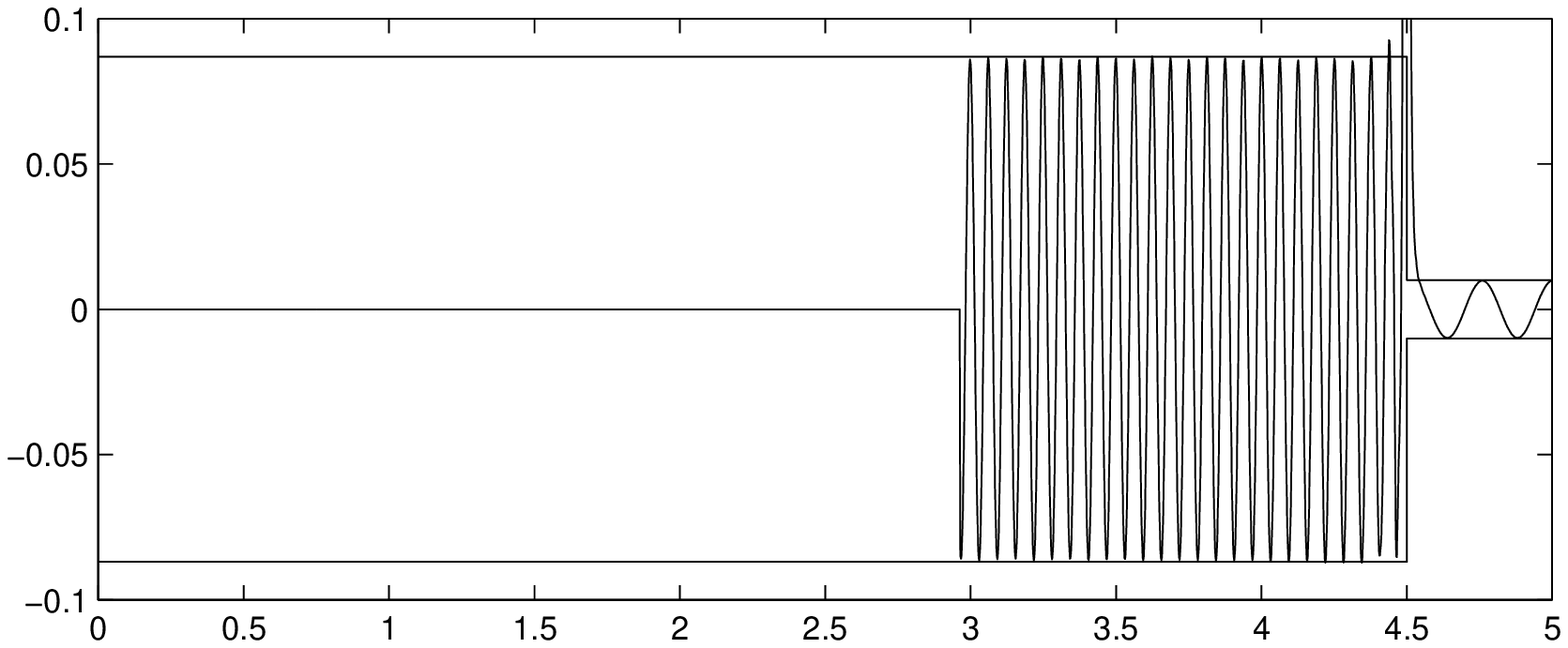}} \\

 \subfigure[$\kappa=39$, N=800, before the interpolation]{
     \label{fig2:subfig:c}
     \includegraphics[width=12cm]{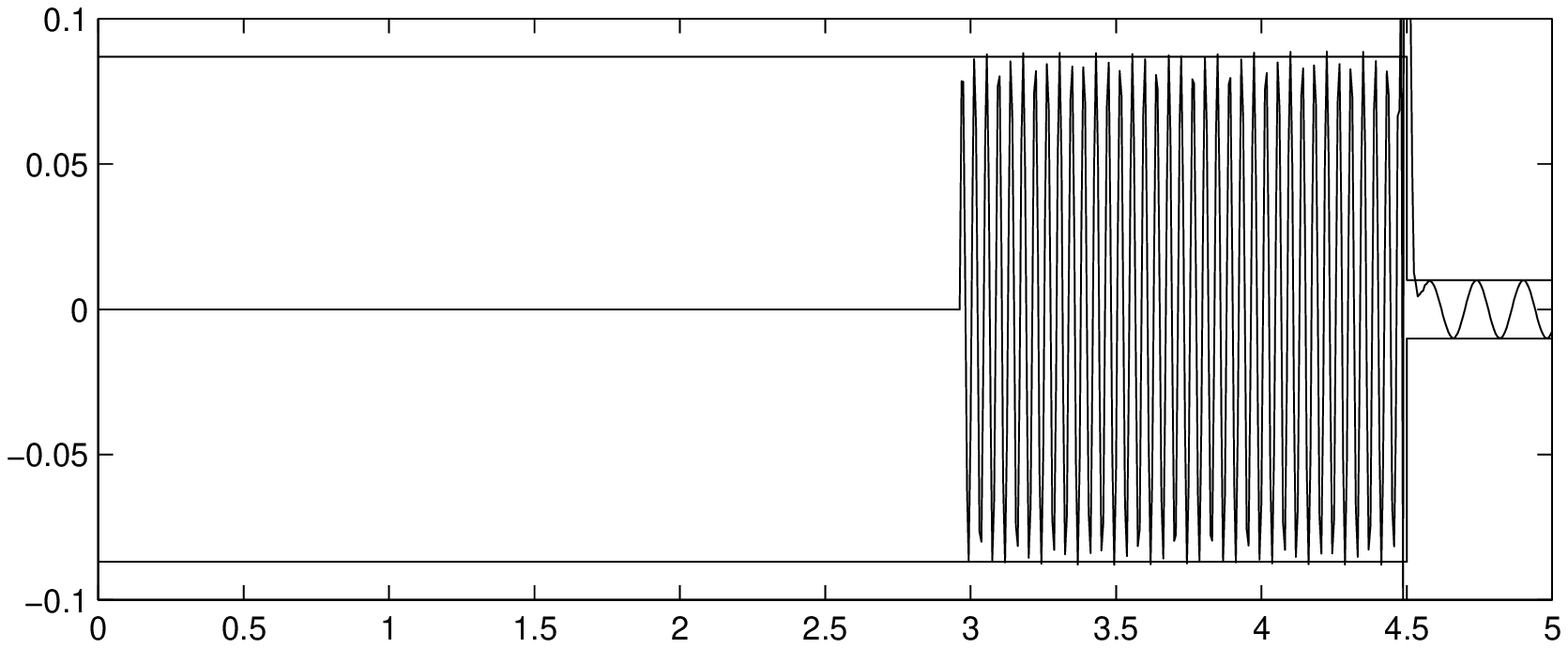}} \\    
      
 \subfigure[$\kappa=39$, N=800, after the interpolation]{
     \label{fig2:subfig:d}
     \includegraphics[width=12cm]{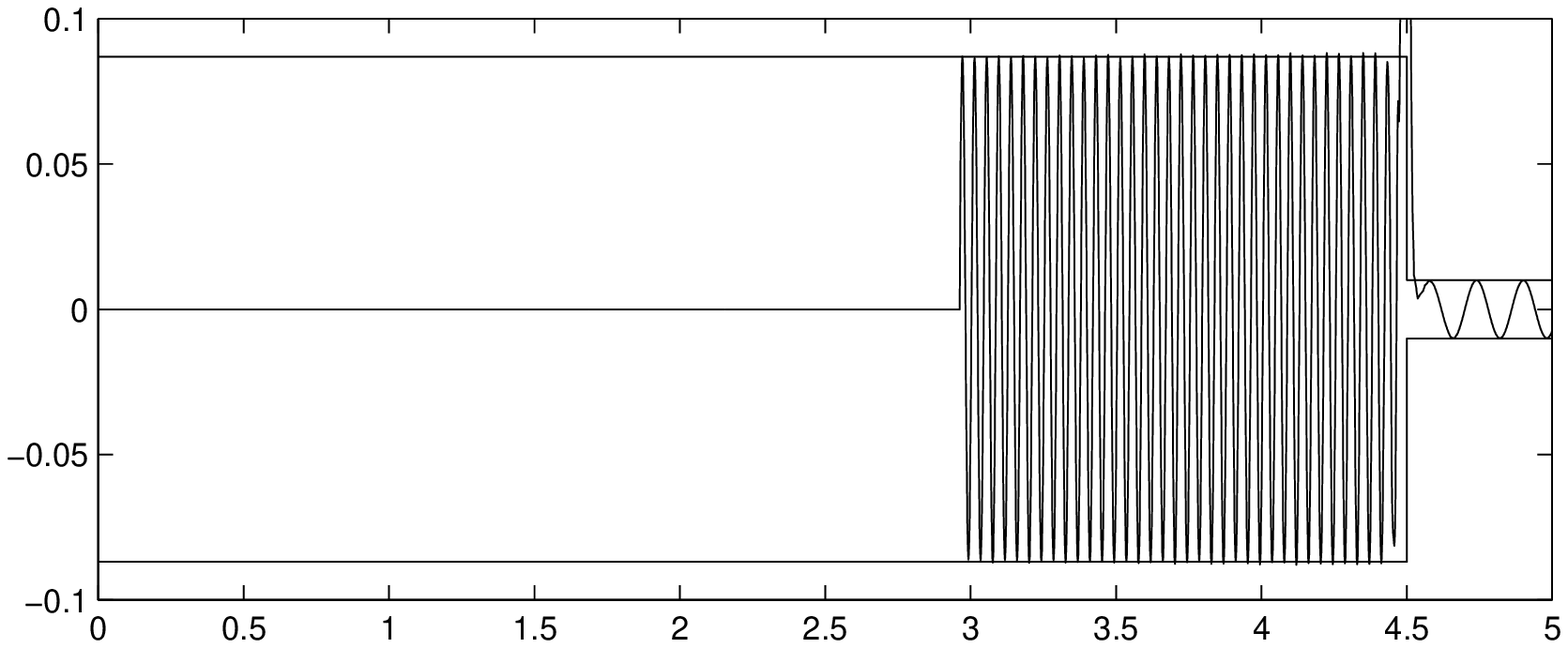}} \\
     
% \subfigure[$\kappa=52$, N=800, before the re-interpolation]{
%     \label{fig2:subfig:e}
%     \includegraphics[width=12cm]{d52b.eps}} \\  
     
% \subfigure[$\kappa=52$, N=800, after the re-interpolation]{
%     \label{fig2:subfig:f}
%     \includegraphics[width=12cm]{d52.eps}} \\  
     
 \subfigure[$\kappa=52$, N=1600]{
     \label{fig2:subfig:g}
     \includegraphics[width=12cm]{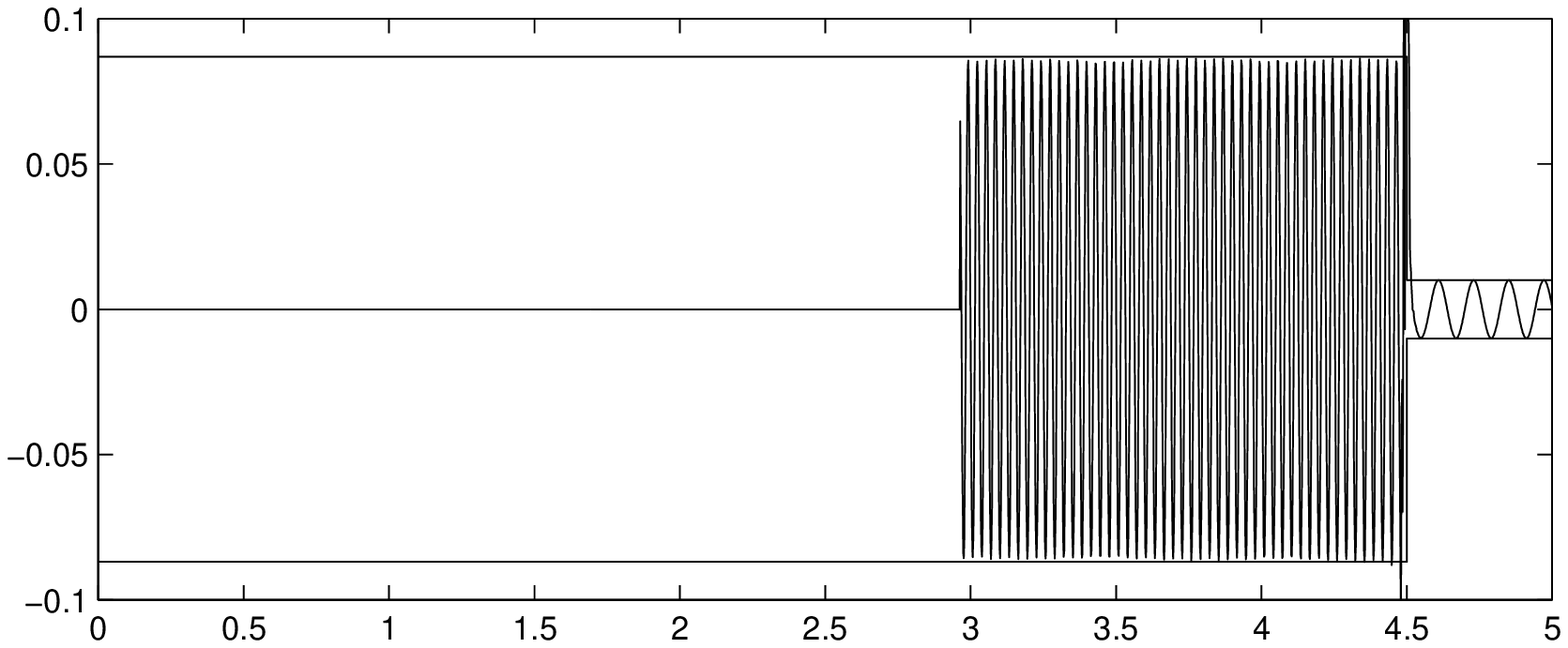}} \\       
     
 \subfigure[$\kappa=65$, N=1600]{
     \label{fig2:subfig:h}
     \includegraphics[width=12cm]{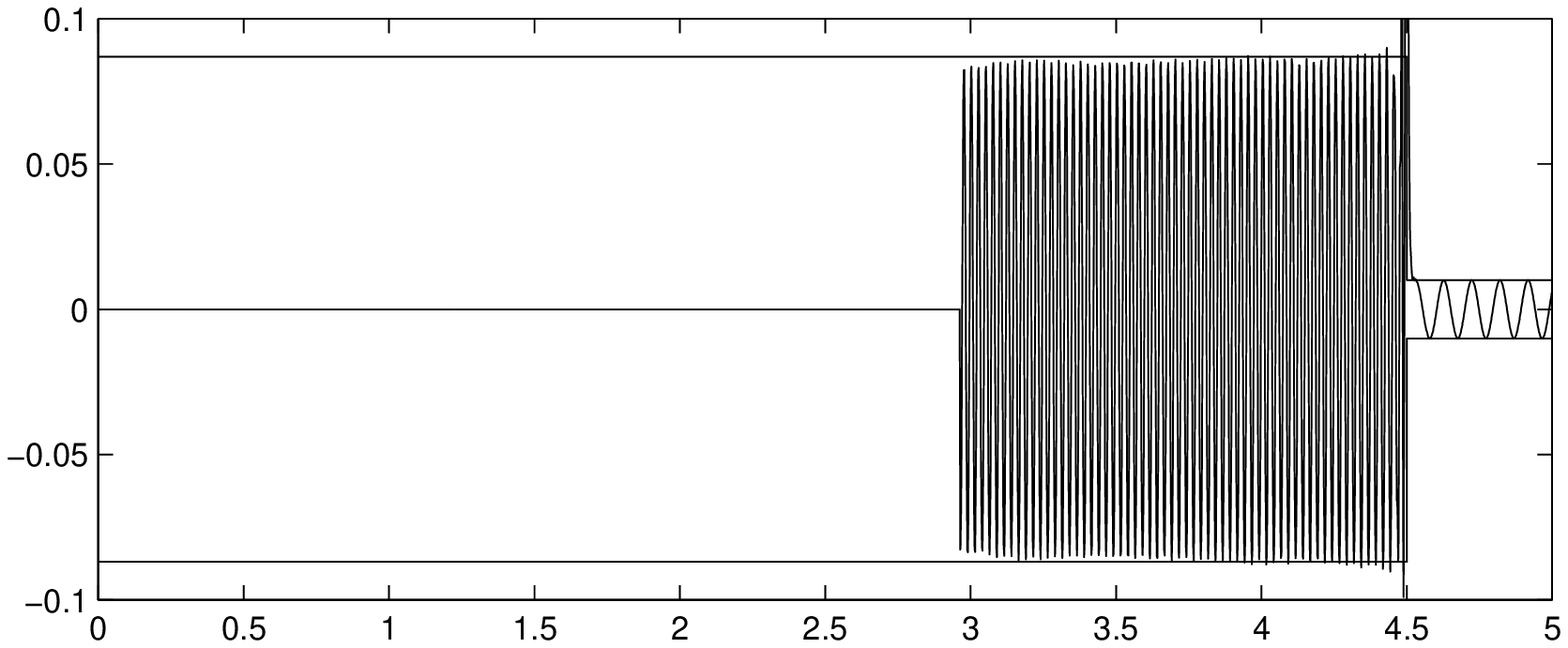}} \\       
     
% \subfigure[$\kappa=78$, N=1600]{
%     \label{fig2:subfig:i}
%     \includegraphics[width=12cm]{d78.eps}} \\       
     
\caption{Shock entropy wave interaction.}
% (a)-(g): Amplitude of the entropy wave} 
\label{shockfig1}
\end{figure}
\end{center}

\end{document}